\documentclass{amsart}
\usepackage{amssymb}
\usepackage{amscd}

\theoremstyle{plain}
\newtheorem{theorem}{Theorem}[section]
\newtheorem{lemma}[theorem]{Lemma}
\newtheorem{proposition}[theorem]{Proposition}
\newtheorem{corollary}[theorem]{Corollary}

\newtheorem*{TA}{Theorem A}
\newtheorem*{TB}{Theorem B}

\theoremstyle{definition}
\newtheorem{definition}[theorem]{Definition}
\newtheorem*{notation}{Notation}
\newtheorem{problem}{Problem}
\newtheorem{remark}[theorem]{Remark}
\newtheorem{example}[theorem]{Example}

\numberwithin{equation}{section}

\newcommand{\dd}{\mathrm{d}}

\newcommand{\hotimes}{\mathbin{\widehat{\otimes}}}
\newcommand{\ootimes}{\mathbin{\ol{\otimes}}}
\newcommand{\jz}{$\set{\jj,0}$}

\newcommand{\bp}{\mathbin{\square}}
\newcommand{\bpd}{\mathbin{\square^\dd}}
\newcommand{\bpz}{\boxdot}

\newcommand{\ltp}{\boxtimes}
\newcommand{\FD}{\mbf{F_D}}
\newcommand{\vx}{\mathsf{x}}
\newcommand{\id}{\mathrm{id}}
\newcommand{\ltr}{{\mathbin{\vartriangleleft}}}
\newcommand{\utr}{{\mathbin{\bigtriangleup}}}
\newcommand{\dtr}{{\mathbin{\bigtriangledown}}}
\newcommand{\Coni}{\Con_{\infty}}

\newcommand{\Free}{\tup{F}}

\newcommand{\fa}{\mathfrak{a}}

\newcommand{\fc}{\mathfrak{c}}
\newcommand{\fg}{\mathfrak{g}}

\newcommand{\ghost}{{{}_{-}}}

\DeclareMathOperator{\Con}{Con}
\DeclareMathOperator{\Conc}{Con_c}

\DeclareMathOperator{\Pow}{\tup{P}}
\DeclareMathOperator{\BoxCl}{Box}

    [1998/03/01 Commands for lattices]
\RequirePackage{amsmath}
\RequirePackage{amssymb}
\RequirePackage{latexsym}
\RequirePackage{eucal}
\RequirePackage{xspace}
\RequirePackage{verbatim}

\newcommand{\jj}{\vee}
\newcommand{\mm}{\wedge}
\newcommand{\JJ}{\bigvee}
\newcommand{\MM}{\bigwedge}
\newcommand{\JJm}[2]{\JJ(\,#1\mid#2\,)}
\newcommand{\MMm}[2]{\MM(\,#1\mid#2\,)}

\newcommand{\uu}{\cup}
\newcommand{\ii}{\cap}
\newcommand{\UU}{\bigcup}
\newcommand{\II}{\bigcap}
\newcommand{\UUm}[2]{\UU(\,#1\mid#2\,)}
\newcommand{\IIm}[2]{\II(\,#1\mid#2\,)}

\newcommand{\ci}{\subseteq}
\newcommand{\sci}{\subset}
\newcommand{\nin}{\notin}
\newcommand{\es}{\varnothing}
\newcommand{\set}[1]{\{#1\}}
\newcommand{\setm}[2]{\{\,#1\mid#2\,\}}
\def\vv<#1>{\langle#1\rangle}

\newcommand{\ga}{\alpha}
\newcommand{\gl}{\lambda}
\newcommand{\gm}{\mu}
\newcommand{\gr}{\varrho}
\newcommand{\gs}{\sigma}

\newcommand{\gF}{\Phi}

\newcommand{\gQ}{\Theta}
\newcommand{\gY}{\Psi}

\newcommand{\tbf}{\textbf}
\newcommand{\trm}{\textrm}
\newcommand{\tup}{\textup}
\newcommand{\mbf}{\mathbf}

\newcommand{\E}[1]{\mathcal{#1}}
\newcommand{\F}[1]{\mathfrak{#1}}

\newcommand{\q}{\quad}
\newcommand{\qq}{\qquad}
\newcommand{\iso}{\cong}
\newcommand{\ol}[1]{\overline{#1}\xspace}

\begin{document}

\title[Box products of lattices]
{A new lattice construction: the box product}

\author{G.~Gr\"atzer}
\thanks{The research of the first author was partially
         supported by the NSERC of Canada.}
\address{Department of Mathematics\\
          University of Manitoba\\
          Winnipeg, Manitoba\\
          Canada R3T 2N2}
\email{gratzer@cc.umanitoba.ca}
  \urladdr{http://www.maths.umanitoba.ca/homepages/gratzer/}

  \author{F.~Wehrung}
  \address{C.N.R.S.\\
           Universit\'e de Caen, Campus II\\
           D\'epartement de Math\'ematiques\\
           B.P. 5186\\
           14032 Caen Cedex\\
           France}
  \email{wehrung@math.unicaen.fr}
  \urladdr{http://www.math.unicaen.fr/\~{}wehrung}
\date{February 16, 1999}

\begin{abstract}
In a recent paper, the authors have proved that
for lattices $A$ and $B$ with zero, the isomorphism
  \[
    \Conc(A \otimes B) \iso \Conc A \otimes\Conc B,
  \]
holds, provided that the tensor product satisfies a very natural
condition (of being \emph{capped}) implying that $A \otimes
B$ is a lattice. In general, $A \otimes B$ is not a lattice;
for instance, we proved that $M_3 \otimes \Free(3)$ is not a lattice.

In this paper, we introduce a new lattice construction, the \emph{box
product} for arbitrary lattices. The tensor product construction for
complete lattices introduced by G.~N. Raney in 1960 and by
R.~Wille in 1985 and the tensor product construction of A.~Fraser in 1978
for semilattices bear some formal resemblance to the new construction.

For lattices $A$ and $B$, while their tensor product $A \otimes B$ (as
semilattices) is not always a lattice, the box product, $A\bp B$, is always a
lattice. Furthermore, the box product and some of its ideals behave like an
improved tensor product. For example, if $A$ and $B$ are lattices with unit,
then the isomorphism
  \[
    \Conc(A \bp B) \iso \Conc A \otimes\Conc B
  \]
  holds. There are analogous results for lattices $A$ and $B$ with zero
and for a bounded lattice $A$ and an arbitrary lattice $B$.

A join-semilattice $S$ with zero is called \emph{$\set{0}$-representable},
if there exists a lattice $L$ with zero such that $\Conc L \iso S$.  The
above isomorphism results yield the following consequence: \emph{The tensor
product of two $\set{0}$-representable semilattices is
$\set{0}$-representable}.
  \end{abstract}

\maketitle

\section{Introduction}
In our paper \cite{GrWe2}, we recalled in detail the introduction of
tensor products of lattices in the seventies.  The main result of this
field is the isomorphism
\begin{equation}\label{E:old}
    \Conc(A \otimes B) \iso \Conc A \otimes\Conc B
\end{equation}
  we proved in \cite{GrWe2} for capped tensor products; this generalizes the
result of G. Gr\"atzer, H. Lakser, and R. W. Quackenbush \cite{GLQu81} for
finite lattices.  This isomorphism does not always make sense because
$A \otimes B$ is not a lattice, in general; in \cite{GrWe3} and \cite{GrWe4},
we provided examples, for instance, $M_3 \otimes \tup{F}(3)$ is not a
lattice (this solved a problem proposed in R.~W. Quackenbush \cite{rQ85}).

In \cite{GrWe4}, we solved a problem of E.~T. Schmidt and the first
author: does every lattice have a proper congruence-preserving extension.
In earlier papers, such an extension for a distributive lattice was
provided by Schmidt's $M_3[D]$ construction.  Trying to use this
construction in the general case ran into the same type of problem
mentioned in the previous paragraph: for a general lattice $L$, the
construction $M_3[L]$ does not always yield a lattice.  The problem was
solved by the $M_3\langle L\rangle$ construction that inherits some
properties of the $M_3[L]$ construction and always produces a lattice.

In this paper, we introduce the \emph{box product} of lattices
(Definition~\ref{D:BoxProd}). For~lattices $A$ and $B$, the box product,
$A \bp B$, is always a lattice.
If $A$ and $B$ are finite, then $A\bp B$ is isomorphic to the
complete tensor product $A\hotimes B$ considered in R.~Wille \cite{Wille},
see also Section~\ref{S:discussion}.

We also introduce an ideal $A \ltp B$ of $A \bp B$; we shall
call $A \ltp B$ the \emph{lattice tensor product} of $A$ and $B$. The ideal
$A\ltp B$ can be defined if $A$ and $B$ have a zero,
or if either $A$ or $B$ is bounded,
or if $A$ and $B$ have unit, see Lemma~\ref{L:ltpexists}. At the
end of Section~\ref{S:stp}, we point out that the lattice tensor product
$M_3 \ltp L$ and $M_3\langle L\rangle$ are isomorphic, showing how the
concept of lattice tensor product was inspired by the $M_3\langle
L\rangle$ construction.

This paper makes the first few steps in exploring the connections among
$A \otimes B$, $A \bp B$, and $A \ltp B$. If $A$ or $B$ is distributive, then
$A\ltp B=A\otimes B$  (Proposition~\ref{P:DLatTens}). The $A\ltp B$
construction yields a universal object for a certain kind of
``bimorphism", see Definition~\ref{D:bimorphism} and
Proposition~\ref{P:TensUniv}. The lattice $A \ltp B$ is always a capped
sub-tensor product of $A$ and $B$ (in the sense of \cite{GrWe2}), see
Theorem~\ref{T:ltpTensLatt}. By using the isomorphism result of
\cite{GrWe2} (see \eqref{E:old}), this yields the isomorphism
  \begin{equation}\label{Eq:ConIsom}
    \Conc(A \ltp B) \iso \Conc A \otimes \Conc B,
  \end{equation}
  see Theorem~\ref{T:AotBCon}. A direct limit argument extends
this isomorphism to two arbitrary lattices, one of which is bounded
(Theorem~\ref{T:AotBCon2}). Finally, as a ``dual'' of
\eqref{E:old}, we prove that if $A$ and $B$ are lattices with unit, then
the isomorphism
  \begin{equation}\label{Eq:unit}
    \Conc(A \bp B) \iso \Conc A \otimes \Conc B
  \end{equation}
holds (Theorem~\ref{T:AotBCon3}).

These isomorphism statements have some interesting consequences related to
the classical Congruence Lattice Characterization Problem; we refer the
reader to \cite{GrSc} for a review of this field. Let us say
that a join semilattice $S$ with zero is \emph{representable}  (resp.,
$\set{0}$-representable, \emph{$\set{0, 1}$-representable}), if there
exists a lattice $L$ (resp., a lattice with zero, a bounded lattice) such
that the join semilattice $\Conc L$ of compact congruences of $L$ is
isomorphic to $S$. In this paper, we prove two related results:

\begin{TA}
Let $S$ and $T$ be $\set{0}$-representable join semilattices. Then
the tensor product $S\otimes T$ is also $\set{0}$-representable.
\end{TA}

\begin{TB}
Let $S$ and $T$ be join semilattices. If $S$ is $\set{0, 1}$-representable
and $T$ is representable, then the tensor product $S \otimes T$ is
representable.
  \end{TB}

We will use the notations and terminology of \cite{GrWe2} and
\cite{GrWe3}. For any set $X$, we shall denote by $\tup{P}(X)$ the power
set of $X$, and $\tup{P}^*(X)=\tup{P}(X) - \set{\es,X}$.

If $L$ is a lattice, the statement ``$0_L$ exists" means that $L$ has a
least element, which we shall always denote by $0_L$; and, similarly, for
$1_L$, the largest element of~$L$.

$\E L_0$ denotes the category of all lattices with zero and
$\set{0}$-homomorphisms. Let $L^\dd$ denote the dual of the lattice $L$.

A non-negative integer $n$ will be identified with the set
$\set{0,1,\ldots,n-1}$.  For a positive integer $n$, let $\Pow(n)$
denote the power set of $n$, partially ordered by inclusion.

Let $L$ be a lattice, let $n>0$, and let $a_0$,\ldots, $a_{n-1} \in L$. For
a subset $X$ of $n$, we write
\begin{align*}
  a^{(X)}&=\JJm{a_i}{i\in X},\\
  a_{(X)}&=\MMm{a_i}{i\in X}.
\end{align*}
For $b\in L$, define $a^{(\es)}\jj b = b$, even though $a^{(\es)}$ is
not defined unless $L$ has a zero.

We shall sometimes denote a finite list $x_0$,\ldots, $x_{n-1}$
by $\vec x$. For example, if the $x_i$-s are elements of a lattice $L$ and
if $P$ is a lattice polynomial with $n$ variables, then we shall write
$P(\vec x)$ for $P(x_0,\ldots,x_{n-1})$.

\section{The box product}

In this section, we introduce the box product and establish some of its
basic properties.  Throughout this section, let $A$ and $B$ be lattices.

Now we define box products:

\begin{definition}\label{D:BoxProd}
For all $\vv<a,b> \in A \times B$, define
  \[
    a \bp b = \setm{\vv<x, y> \in A \times B}{x \leq a \text{ or } y \leq
       b}.
  \]
  We define the \emph{box product} of $A$ and $B$, denoted by $A \bp B$, as
the set of all finite \emph{intersections} of the form
  \[
     H = \IIm{a_i \bp b_i}{i < n},
  \]
  where $n$ is a positive integer, and $\vv<a_i, b_i> \in A \times B$, for
all $i < n$.
  \end{definition}

$A \bp B$ is a poset under set containment.

  \begin{remark}\label{R:1AboxB}
  It is easy to see that $A \bp B$ has a unit element, $1_{A\bp B}$,
if and only if either $A$ or $B$ does. For example, if $A$ has a unit,
$1_A$, then $1_{A \bp B} = 1_A \bp b$, for all $b \in B$.
  \end{remark}

It is obvious that $A \bp B$ is a meet-subsemilattice of the powerset
lattice of $A \times B$. We shall show in Proposition~\ref{P:bplattice} that
$A \bp B$ is a lattice. First, we need another definition:

\begin{definition}
For $\vv<a, b> \in A \times B$, define
  \[
    a \circ b = \setm{\vv<x, y> \in A \times B} {x \leq a \text{ and }y
    \leq b}.
  \]
We define $A\bpz B$ to be the set of all finite \emph{unions}
of the form
  \begin{equation}\label{Eq:basicH}
    H = \UUm{a_i \bp b_i}{i < m} \uu \UUm{c_j \circ d_j}{j < n},
  \end{equation}
where $m > 0$ and $n \geq 0$ are integers, $a_i$, $c_j \in A$, and
$b_i$, $d_j \in B$.
  \end{definition}

The proof of the following lemma is straightforward; the details are left
to the reader.

\begin{lemma}\label{L:BoxCircId}
Let $a$, $a' \in A$ and $b$, $b' \in B$. Then the following
assertions hold:
  \begin{itemize}
  \item[\trm{(a)}] $a \circ b \ci a' \bp b'$ if and only if
  $a \leq a'$ or $b \leq b'$.

  \item[\trm{(b)}]
  $(a \circ b) \ii (a' \circ b') = (a \mm a') \circ (b \mm b')$.

  \item[\trm{(c)}]
  $(a \bp b) \ii (a' \circ b') =
  ((a \mm a') \circ b') \uu (a' \circ (b \mm b'))$.

  \item[\trm{(d)}]
  $(a \bp b) \ii (a' \bp b') =
  ((a \mm a') \bp (b \mm b')) \uu (a \circ b') \uu (a' \circ b)$.

\item[\trm{(e)}] $a \bp b \ci a' \bp b'$ if and only if either
$A = (a']$, or $B = (b']$, or ($a \leq a'$ and $b \leq b'$).
  \end{itemize}
  \end{lemma}

\begin{corollary}\label{C:AbpzBsubl}
$A \bpz B$ is a sublattice of $\Pow(A \times B)$.
  \end{corollary}

Let $L$ be a lattice; a \emph{closure system} on $L$ is a
subset $K$ of $L$ such that for every element $x$ of $L$,
there exists a least element $\ol{x}$ of $K$ satisfying $x \leq \ol{x}$.
Note that $K$ is then automatically a meet-subsemilattice of $L$.
The element $\ol{x}$ is called the \emph{closure of $x$ in $K$}.

The following well-known lemma requires no proof.

\begin{lemma}\label{L:ClosSys}
Let $L$ be a lattice and let $K$ be a closure system on $L$.
Then $K$ is a lattice and the join is given by the formula
  \begin{equation*}
  x\jj_K y = \ol{x\jj_L y}.
  \end{equation*}
  \end{lemma}

The following lemma is fundamental in the theory of box products.

\begin{lemma}\label{L:BoxClos}
$A \bp B$ is a closure system in $A \bpz B$.
  \end{lemma}

\begin{proof}
  Let
  \begin{equation}\label{Eq:DecH}
    H =  \UUm{a_i \bp b_i}{i<m} \uu \UUm{c_j \circ d_j}{j<n} \in A \bpz B,
  \end{equation}
where $m>0$ and $n \geq 0$. Put $\ol{a} = \JJm{a_i}{i<m}$ and
$\ol{b} = \JJm{b_i}{i < m}$. Set
  \[
    \ol{H} = \IIm{(\ol{a}\jj c^{(X)})\bp(\ol{b}\jj d^{(n-X)})}{X \ci n}.
  \]
  Note that $\ol{H} \in A \bp B$. We shall prove that $\ol{H}$ is the
closure of $H$ in $A \bp B$.

First, we verify that $H \ci \ol{H}$. For all $i<m$, $a_i \bp b_i \ci\ol{a}
\bp\ol{b} \ci \ol{H}$.

Let $j<n$ and let $X \ci n$; we prove that
$c_j \circ d_j \ci c^{(X)}\bp d^{(n-X)}$.
If $j \in X$, then $c_j \leq c^{(X)}$, and so the conclusion
follows by Lemma~\ref{L:BoxCircId} (a). Similarly, if $j\notin X$, then
$d_j \leq d^{(n - X)}$, and so the conclusion follows again by
Lemma~\ref{L:BoxCircId} (a). In both cases, $c_j \circ d_j \ci \ol{H}$.
Hence $H \ci \ol{H}$.

Second, it suffices to prove that for all $\vv<a,b> \in A \times B$,
$H \ci a \bp b$ implies that $\ol{H} \ci a \bp b$. This conclusion is
trivial if $A = (a]$ or if $B = (b]$, so suppose that $a$ (resp., $b$) is
not the greatest element of $A$ (resp., of $B$). For all $i<m$, the
containment $a_i \bp b_i \ci H \ci a \bp b$ holds, thus, by
Lemma~\ref{L:BoxCircId} (e), $a_i \leq a$ and $b_i \leq b$; it follows
that $\ol{a} \leq a$ and $\ol{b} \leq b$. Put $X = \setm{j<n}{c_j \leq
a}$. Since $\ol{a} \leq a$, it follows from the definition that
$c^{(X)} \leq a$. Furthermore,  $c_j\nleq a$, for all $j \in n - X$; but
$\vv<c_j,d_j> \in H \ci a \bp b$, thus $d_j \leq b$. It follows that
$d^{(n -  X)} \leq b$. Therefore,
$\ol{H} \ci (\ol{a}\jj c^{(X)})\bp(\ol{b}\jj d^{(n-X)}) \ci a \bp b$.
  \end{proof}

We shall call $\ol{H}$ the \emph{box closure} of $H$ and denote it by
$\BoxCl(H)$. Since $\BoxCl(H)$ is the least element of $A\bp B$ containing
$H$, it is independent of the decomposition \eqref{Eq:DecH}. This
definition can be extended to all subsets of $A\times B$:

\begin{definition}\label{D:BoxClos}
Let $A$ and $B$ be lattices. For $X\ci A\times B$, we define the
\emph{box closure} of $X$:
  \[
  \BoxCl(X)=\IIm{a\bp b}{\vv<a,b>\in A\times B,\ X\ci a\bp b}.
  \]
\end{definition}

So the box closure of $X$ is the intersection of all elements of $A\bp B$
containing~$X$.  For an arbitrary subset $X$ of $A\times B$, it may not
belong to $A\bp B$.

  \begin{proposition}\label{P:bplattice}
   Let $A$ and $B$ be lattices. If $H \in A \bpz B$, then
  $\BoxCl(H) \in A \bp B$. In particular, $A\bp B$ is a lattice.
  \end{proposition}

It is important to note that the proof of Lemma~\ref{L:BoxClos} gives us the
existence of $\BoxCl(H)$, for $H\in A\bpz B$, as well as effective
\emph{formulas} to compute $\BoxCl(H)$.

The following definition is motivated by R.~Wille \cite{Wille}:

\begin{definition}\label{D:utretc}
Let $A$ and $B$ be lattices.
\begin{enumerate}
\item For $a$, $a'\in A$ and $b$, $b'\in B$, we define
  \[
  \vv<a,b>\ltr\vv<a',b'>,\quad\text{if }a\leq a'\text{ or }b\leq b'.
  \]

\item For a subset $X$ of $A\times B$, we define
  \begin{align*}
  X^\utr&=
  \setm{\vv<a,b>\in A\times B}
  {\vv<x,y>\ltr\vv<a,b>,\text{ for all }\vv<x,y>\in X},\\
   X^\dtr&=
  \setm{\vv<a,b>\in A\times B}
  {\vv<a,b>\ltr\vv<x,y>,\text{ for all }\vv<x,y>\in X}.
  \end{align*}
\end{enumerate}
\end{definition}

In particular, $\vv<a,b>\ltr\vv<a',b'>$ if{}f $\vv<a,b>\in a'\bp b'$.
It is easy to characterize the box product and the box closure
in terms of the $\ltr$ relation:

\begin{proposition}\label{P:BoxProdLtr}
Let $A$ and $B$ be lattices. Then
  \[
  A\bp B=\setm{X^\dtr}{X\ci A\times B,\ X\text{ finite}}.
  \]
Furthermore, $\BoxCl(X)=(X^\utr)^\dtr$, for all $X\ci A\times B$.
\end{proposition}

Note the following trivial corollary of Lemma~\ref{L:BoxCircId}(d):

\begin{proposition}\label{P:partbox1}
Every element of $A\bp B$ contains a pure box.
\end{proposition}

The formulas given in Lemma~\ref{L:BoxClos} to compute the box closure of
an element of $A \bpz B$ can be used to give direct expressions for the
join of two elements of $A \bp B$, as follows. For all positive integers
$m$ and $n$, let $\gs_{m, n}$ be an effectively constructed bijection from
$2^m + 2^n - 4$ onto the ``disjoint union'' of $\tup{P}^*(m)$ and
$\tup{P}^*(n)$, that is, onto
$(\tup{P}^*(m) \times \set{0}) \uu (\tup{P}^*(n) \times \set{1})$.
For all $k < 2^m + 2^n - 4$, we define the
lattice polynomials $M_{m, n, k}$ and $N_{m, n, k}$ by
  \begin{align}
    M_{m, n, k}(\vec a,\vec c)& =
  \begin{cases}
    \MMm{a_i}{i \in X},& \text{ if }\gs_{m, n}(k) = \vv<X,0>;\\
    \MMm{c_j}{j \in Y},& \text{ if }\gs_{m, n}(k) = \vv<Y,1>;
\end{cases}
\label{Eq:Mmnk}
  \intertext{and}
N_{m, n, k}(\vec b,\vec d)& =
  \begin{cases}
    \MMm{b_i}{i \in m -  X},& \text{ if }\gs_{m, n}(k) = \vv<X,0>;\\
    \MMm{d_j}{j \in n -  Y},& \text{ if }\gs_{m, n}(k) = \vv<Y,1>.
  \end{cases}
\label{Eq:Nmnk}
  \end{align}
Furthermore, for all $\es \ci Z \ci 2^m + 2^n - 4$, we define the lattice
polynomials $U_{m, n, Z}$ and $V_{m, n, Z}$ by the following formulas:
  \begin{align}
    U_{m, n, Z}(\vec a,\vec c) = \MM_{i<m}a_i\jj\MM_{j<n}c_j\jj
  \JJ_{k \in Z}M_{m, n, k}(\vec a,\vec c),\label{Eq:UmnZ}\\
  \intertext{and}
  V_{m, n, Z}(\vec b,\vec d) = \MM_{i<m}b_i\jj\MM_{j<n}d_j\jj
  \JJ_{k \nin Z}N_{m, n, k}(\vec b,\vec d).\label{Eq:VmnZ}
  \end{align}
  By definition, for the cases $Z = \es$ and $Z = 2^m + 2^n - 4$, these
formulas mean:
  \begin{align}
    U_{m, n, \es}(\vec a,\vec c)& = \MM_{i<m}a_i\jj\MM_{j<n}c_j,
        \label{Eq:UmnZ0}\\
V_{m, n, 2^m + 2^n - 4}(\vec b,\vec d)& = \MM_{i<m}b_i\jj\MM_{j<n}d_j.
\label{Eq:VmnZ0}
  \end{align}

Now we formulate how the join in $A \bp B$ can be computed:

\begin{lemma}\label{L:CompJoin}
  Let $A$ and $B$ be lattices. Let $H$ and $K \in A \bp B$ be written in
the form
  \begin{align*}
    H &= \IIm{a_i \bp b_i}{i<m},\\
    K &= \IIm{c_j \bp d_j}{j<n}.
  \end{align*}
Then
  \[
    H\jj K =
\IIm{U_{m, n, Z}(\vec a,\vec c) \bp V_{m, n, Z}(\vec b,\vec d)}
{Z \ci 2^m+2^n-4}.
  \]
  \end{lemma}

\begin{proof}
A direct computation shows that
\begin{align*}
H \uu K = \left(\MMm{a_i}{i<m} \bp \MMm{b_j}{j<n}\right) \uu
\left(\MMm{c_i}{i<m} \bp\MMm{d_j}{j<n}\right)\\
   \uu \UUm{M_{m, n, k}(\vec a,\vec c) \circ N_{m, n, k}(\vec b,\vec d)}
{k<2^m+2^n-4}.
\end{align*}
The conclusion follows right away from the proof of
Lemma~\ref{L:BoxClos} and the definition of the
polynomials $U_{m, n, Z}$, $V_{m, n, Z}$.
\end{proof}

\section{Pure lattice tensors; lattice tensor product}

\begin{definition}\label{D:BotPureltp}
Let $A$, $B$, and $L$ be lattices.
\begin{enumerate}
\item We define the \emph{bottom} of $L$ by
  \[
  \bot_L=
   \begin{cases}
    \set{0_L},&\text{if }L\text{ has a zero};\\
    \es,&\text{otherwise}.
   \end{cases}
  \]
\item We put
  \[
  \bot_{A,B}=(A\times\bot_B)\uu(\bot_A\times B).
  \]

\item Let $\vv<a,b>\in A\times B$.
We define the \emph{pure lattice tensor} of $a$ and $b$:
  \[
  a\ltp b=(a\circ b)\uu\bot_{A,B}.
  \]

\item A subset $X$ of $A\times B$ is \emph{confined}, if $X\ci a\ltp b$, for
some $\vv<a,b>\in A\times B$.

\item A subset $H$ of $A\times B$ is a \emph{bi-ideal} of $A\times B$, if
the following conditions hold:
\begin{enumerate}
\item $\bot_{A,B}\ci H$;
\item $H$ is a hereditary subset of $A\times B$;
\item For $a_0$, $a_1\in A$ and $b\in B$, if $\vv<a_0,b>\in H$ and
$\vv<a_1,b>\in H$, then $\vv<a_0\jj a_1,b>\in H$; and symmetrically.
\end{enumerate}

\end{enumerate}
\end{definition}

As an immediate consequence of the definition of a bi-ideal, we obtain:

  \begin{lemma}\label{L:biideal}
Let $A$ and $B$ be lattices. The elements of $A \bp B$ are bi-ideals of
$A \times B$.
  \end{lemma}

Now the lattice tensor product:

\begin{definition}\label{D:LattTP}
Let $A$ and $B$ be lattices. Let $A\ltp B$ be the set of all confined
elements of $A\bp B$. If $A\ltp B$ is nonempty, then we say that $A\ltp B$
is \emph{defined}, and we call it the \emph{lattice tensor product} of $A$
and $B$.
\end{definition}

We obtain immediately the following trivial consequence of
Definitions \ref{D:BotPureltp} and \ref{D:LattTP}:

\begin{proposition}\label{P:AltpBlatt}
Let $A$ and $B$ be lattices. If $A\ltp B$ is defined, then it is an
ideal of $A\bp B$. In particular, it is a lattice.
\end{proposition}

  Note that if $A$ and $B$ have zero, then $a \ltp b$ is the same as
$a\otimes b$ in \cite{GrWe2}. However, the underlying structures, $A \ltp B$
(see Definition~\ref{D:LattTP}) and $A\otimes B$ (see \cite{GrWe2}) are
different.

Note the following trivial corollary of Proposition~\ref{P:partbox1}:

\begin{proposition}\label{P:partbox2}
Every element of $A\ltp B$ contains a (confined) pure box.
\end{proposition}

Now we completely characterize when $A \ltp B$ is defined:

\begin{lemma}\label{L:ltpexists}
  Let $A$ and $B$ be lattices. Then $A \ltp B$ is defined
if{}f one of the
following conditions hold:
  \begin{enumerate}
  \item $A$ and $B$ are lattices with zero;
  \item $A$ and $B$ are lattices with unit;
  \item $A$ or $B$ is bounded.
\end{enumerate}
  \end{lemma}

\begin{proof}
  Let (i) hold. Let $\vv<a,b> \in A \times B$. Then
  \[
    a \ltp b=(a \bp 0_B) \ii (0_A \bp b).
  \]
  Therefore, $a \ltp b \in A \bp B$ and it is confined (by itself).
Thus $a \ltp b \in A \ltp B$ and so $A \ltp B \neq \es$.

Let (ii) hold. Then every element of $A \bp B$ is confined by
$1_A \ltp 1_B = A \times B$, and so $A \ltp B = A \bp B \neq \es$.

Let (iii) hold. If $A$ is a bounded lattice and $b \in B$, then
$0_A \bp b = 1_A \ltp b$,
hence $A \ltp B \neq \es$. If $B$ is a bounded
lattice, we proceed symmetrically.

Now, conversely, let us assume that $A \ltp B$ is defined, that is,
$A \ltp B \neq \es$.  There are 16 cases to consider whether $A$ and $B$ have
zero and/or unit.  Nine of these possibilities are covered by (i)--(iii); the
remaining seven possibilities, by symmetry, are covered by the following
single case:

The lattice $A$ has no zero and the lattice $B$ has no unit. If
$H \in A \ltp B$ is confined by $a \ltp b$, $a \in A$, $b \in B$, then there
is a pure box $u\bp v$ confined by $a \ltp b$,
by Proposition~\ref{P:partbox2}. Since $A$ has no
zero, $u \in A^-$. Thus $\vv<u, x> \leq \vv<a, b>$, for all $x \in B$; hence
$b$ is the unit of $B$, a contradiction.
  \end{proof}

Box closures play an important role for lattice tensor products:

\begin{lemma}\label{L:bddclosure}
  Let $A$ and $B$ be lattices.
\begin{enumerate}
  \item For $a \in A$ and $b \in B$,
  \[
    \BoxCl(a \ltp b) = a \ltp b.
  \]
\item Let $H \ci A \times B$. If $H$ is confined, then the box closure of
$H$ is also confined.
  \item $K \in A \ltp B$ if{}f $K$ is the box closure of some confined
$H \in A \bpz B$.
  \item  If $A$ and $B$ are lattices with zero and $a_0$, $a_1 \in A$,
$b_0$, $b_1 \in B$ satisfy $a_0 \leq a_1$ and $b_0 \leq b_1$, then
  \[
    \BoxCl((a_0 \ltp b_1) \uu (a_1 \ltp b_0)) =
    (a_0 \ltp b_1) \uu (a_1 \ltp b_0),
  \]
  so $(a_0 \ltp b_1) \uu (a_1 \ltp b_0) \in A \ltp B$.
\end{enumerate}
  \end{lemma}

\begin{proof}\hfill

(i). Since
  \[
    \IIm{x\bp b}{x\in A}\ii\IIm{a\bp y}{y\in B}=a\ltp b,
  \]
  it follows that $\BoxCl(a\ltp b)=a\ltp b$.

(ii). If $H$ is confined, then there exists $\vv<a,b>\in A\times B$ such
that $H \ci a\ltp b$. Therefore, $\BoxCl(H) \ci \BoxCl(a\ltp b) = a\ltp b$,
by (i). So $\BoxCl(H)$ is confined.

(iii). If $K \in A \ltp B$, then $K \in A \bp B$ and $K$ is confined, so it
is the box closure of some confined $H \in A \bpz B$, namely, of $H = K$.
Conversely, the box closure $K$ of any confined $H \in A \bpz B$ is in
$A \bp B$ and, by (i), it is confined, hence $K \in A \ltp B$.

(iv). This follows from the formula:
  \[
    (a_0 \ltp b_1) \uu (a_1 \ltp b_0) = (a_0 \bp b_0) \ii (0_A \bp b_1) \ii
(a_1 \bp 0_B).
  \]

\vspace{-15pt}
  \end{proof}

For lattices $A$ and $B$ with unit, every subset of $A\times B$ is confined
(by $1_A\ltp 1_B=A\times B$). In particular, $A\ltp B=A\bp B$. For the two
other cases of Lemma~\ref{L:ltpexists}, we describe the elements of
$A\ltp B$:

\begin{lemma}\label{L:withzero}
  Let $A$ and $B$ be lattices with zero. Then the elements of $A\ltp B$ are
exactly the finite intersections of the form
  \begin{align}\label{Eq:ReprM0}
  H=\IIm{a_i\bp b_i&}{i<n},\\
\intertext{satisfying}
  \MMm{a_i}{i<n}&=0_A,\notag\\
   \MMm{b_i}{i<n}&=0_B,\notag
  \end{align}
where $n>0$, $\vv<a_i,b_i>\in A\times B$, for all $i<n$. Furthermore, every
element of $A\ltp B$ can be written as a finite \emph{union} of pure lattice
tensors:
  \begin{equation}\label{Eq:uaibi}
  H=\UUm{a_i \ltp b_i}{i<n},
  \end{equation}
where $x\in B$, $n\geq0$, and $\vv<a_i,b_i>\in A\times B$, for all $i<n$.

Conversely, the box closure of any element of the form \eqref{Eq:uaibi}
belongs to $A\ltp B$.
  \end{lemma}

It follows, in particular, that the elements of $A\ltp B$ are exactly the
elements of the form $\JJm{a_i\ltp b_i}{i < n}$, where $n>0$, $a_0$,\dots,
$a_{n-1}\in A$, and $b_0$,\dots, $b_{n-1}\in B$, that is, the pure
lattice tensors form a join-basis of $A\ltp B$.

\begin{proof}
Let $H\in A\bp B$. If $H\in A\ltp B$, then there exists
$\vv<a,b>\in A\times B$ such that $H\ci a\ltp b$. Since
$a\ltp b = (0_A\bp b)\ii(a\bp 0_B)$, it follows that
  \[
    H=H\ii(0_A\bp b)\ii(a\bp 0_B)
  \]
can be expressed in the form \eqref{Eq:ReprM0}. Conversely, assume that $H$
is of the form \eqref{Eq:ReprM0}. Observe that
  \[
  a_i\bp b_i=((a_i]_A\times B)\uu(A\times(b_i]_B),
  \]
for all $i<n$. Using the notations $a_{(X)}$ and $b_{(X)}$ (see the
Introduction), we obtain that
  \begin{equation}\label{E:decomp}
    H=(a_{(n)}\bp b_{(n)})\uu\UUm{a_{(X)}\circ b_{(n-X)}}
     {\es\sci X\sci n}.
  \end{equation}
  By assumption, $a_{(n)}=0_A$ and $b_{(n)}=0_B$, so we have obtained $H$
as in \eqref{Eq:uaibi}.

Finally, if $H$ is of the form \eqref{Eq:uaibi}, then
$H = (0_A \bp 0_B) \uu \UUm{a_i \circ b_i}{i<n}$, so $H \in A\bpz B$; thus
$\BoxCl(H) \in A\bp B$, by Proposition~\ref{P:bplattice}. Since $H$
is confined (by $u\ltp v$, where $u=\JJm{a_i}{i<n}$ and
$v=\JJm{b_i}{i<n}$), $\BoxCl(H)$ is confined, by Lemma~\ref{L:bddclosure}.
Hence, $\BoxCl(H)$ belongs to $A\ltp B$.
  \end{proof}

The analogue of Lemma~\ref{L:withzero} for the case where $A$ is bounded is
the following:

\begin{lemma}\label{L:Abdd}
  Let $A$ and $B$ be lattices. If $A$ is bounded, then the elements of
$A\ltp B$ are exactly the finite intersections of the form
  \begin{align}\label{Eq:ReprM1}
     H=\IIm{a_i\bp b_i}{i<n},\\
  \intertext{subject to the condition}
    \MMm{a_i}{i<n}=0_A,\notag
  \end{align}
where $n>0$, $\vv<a_i,b_i>\in A\times B$ for all $i<n$. Furthermore, every
element of $A\ltp B$ can be written as a finite \emph{union}
  \begin{align}\label{Eq:0xuuaibi}
    H&=(0_A \bp x) \uu \UUm{a_i \ltp b_i}{i<n}\\
     &=(0_A \bp x) \uu \UUm{a_i \circ b_i}{i<n},\notag
  \end{align}
where $x\in B$, $n\geq0$, and $\vv<a_i,b_i>\in A\times B$, for all $i<n$.

Conversely, the box closure of any element of the form \eqref{Eq:0xuuaibi}
belongs to $A \ltp B$.  The box closures of elements of the form $(0_A \bp x)
\uu (a \ltp b)$ form a join-basis of $A \ltp B$.

  \end{lemma}

\begin{proof}
Let $H\in A\bp B$. If $H\in A\ltp B$, then there exists
$\vv<a,b>\in A\times B$ such that $H\ci a\ltp b$. Since
$a\ltp b\ci 0_A\bp b$,
  \[
    H=H\ii(0_A\bp b)
  \]
can be expressed in the form \eqref{Eq:ReprM1}. Conversely, assume that
$H$ is of the form \eqref{Eq:ReprM1}. Now we proceed as in the proof of
Lemma~\ref{L:withzero} and obtain \eqref{E:decomp}. By assumption,
$a_{(n)}=0_A$, so we are done.

By Proposition~\ref{P:bplattice}, the box closure $\BoxCl(H)$ of any
element $H$ of $A\bpz B$  belongs to $A\bp B$.
Any $H$ of the form \eqref{Eq:0xuuaibi} belongs to $A\bpz B$,
hence $\BoxCl(H)\in A\bp B$. Since $H$ is confined (by $1\ltp v$, where
$v=x\jj\JJm{b_i}{i<n}$), it follows that $\BoxCl(H)$ is confined, by
Lemma~\ref{L:bddclosure}. Hence, $\BoxCl(H)$ belongs to $A\ltp B$.
  \end{proof}

\section{The tensor product of lattices with zero}

By Lemma~\ref{L:withzero}, if $A$ and $B$ are lattices with zero, then
$A\ltp B$ is the set of all box closures of finite subsets of $A\times B$.
Therefore, by Proposition~\ref{P:BoxProdLtr}, we deduce the following:

\begin{proposition}\label{P:ltp0ltr}
Let $A$ and $B$ be lattices with zero. Then
  \[
  A\ltp B=\setm{(X^\utr)^\dtr}{X\ci A\times B,\ X\text{ finite}}.
  \]
\end{proposition}

\begin{corollary}\label{C:BpLtpDual}
Let $A$ and $B$ be lattices with zero. Then
  \[
  A^\dd\bp B^\dd\iso(A\ltp B)^\dd.
  \]
\end{corollary}

\begin{proof}
The pair of maps $X\mapsto X^\utr$, $X\mapsto X^\dtr$ defines a Galois
correspondence between subsets of $A\times B$ (associated with the binary
relation $\ltr$, see Definition~\ref{D:utretc}). Therefore, the second map
defines an isomorphism from the structure
  \begin{align*}
  A^\dd\bp B^\dd&=\setm{X^\utr}{X\ci A\times B,\ X\text{ finite}},
  &&\text{ endowed with containment}\\
  \intertext{onto the structure}
  A\ltp B&=\setm{(X^\utr)^\dtr}{X\ci A\times B,\ X\text{ finite}},
  &&\text{ endowed with reverse containment}.
  \end{align*}
This observation concludes the proof.
\end{proof}

\begin{remark}\label{R:BpLtpDual}
It is easy to describe explicitly the isomorphism in
Corollary~\ref{C:BpLtpDual}. For $\vv<a,b>\in A\times B$, let $a\bpd b$ be
the pure box of $a$ and $b$ in the lattice $A^\dd\bp B^\dd$. Note that
$a\bpd b=\set{\vv<a,b>}^\utr$. Thus, the image of $a\bpd b$ under the
isomorphism of Corollary~\ref{C:BpLtpDual} is
$(\set{\vv<a,b>}^\utr)^\dtr=a\ltp b$. More generally, for a positive integer
$n$ and elements for $\vv<a_i,b_i>\in A\times B$, where $i<n$, the image of
the element
  \[
  \IIm{a_i\bpd b_i}{i<n}
  \]
is
  \[
  \JJm{a_i\ltp b_i}{i<n}
  \]
(the join is computed in $A\ltp B$).
\end{remark}

\begin{notation}
An \emph{upper subset} of a poset $P$ is a subset $X$ with the
property that if $p \in X$ and $p \leq q$ in $P$, then $q \in X$. Let
$\FD(n)$ be the set of all upper subsets
$\fa$ of $\Pow(n)$ such that $\es\notin\fa$ and $n \in\fa$.
\end{notation}

For every upper subset $\fa$ of $\Pow(n)$, note that $\es \nin \fa$
means that $\fa \ne \Pow(n)$, while $n \in\fa$ means that $\fa \ne \es$.
It is easy to see that $\FD(n)$ is a lattice, a sublattice of the
power set of $\Pow(n)$; it is the free distributive lattice
on $n$ generators, where the $i^\mathrm{th}$ generator corresponds to the
element
  \[
    \fg_i = \setm{X  \in \Pow(n)}{i  \in  X}.
  \]

\begin{notation}
For every positive integer $n$ and every $\fa \in \FD(n)$, define
  \[
    \fa^* = \setm{X  \in \Pow(n)}{n -  X \nin \fa}.
  \]
This is similar to the notation $\F{T}^{\#}$ used in Section~3 of R.~Wille
\cite{Wille}.

Furthermore, we associate with $\fa \in \FD(n)$ a lattice polynomial
$P_\fa$, defined by the formula
  \[
    P_\fa(\vx_0,\ldots,\vx_{n-1}) =
     \MMm{\JJm{\vx_i}{i  \in  X}}{X  \in \fa}.
  \]
  \end{notation}

\begin{lemma}\label{L:RepTens}
Let $A$ and $B$ be lattices with zero.
Let $n$ be a positive integer, let $a_0$, \dots, $a_{n-1}$ be elements of
$A$, and let $b_0$, \dots, $b_{n-1}$ be elements of $B$. Then the box
closure of the element
  \[
    H = \UUm{a_i \ltp b_i}{i<n}
  \]
is given by the formula
  \begin{equation}\label{Eq:clos}
    \BoxCl(H) =
     \UUm{P_\fc(a_0,\ldots,a_{n-1}) \ltp
     P_{\fc^*}(b_0,\ldots,b_{n-1})}{\fc  \in \FD(n)}.
  \end{equation}
  \end{lemma}

\begin{proof}
The formulas given in Lemma~\ref{L:BoxClos} for computing
$\BoxCl(H)$ easily give the box closure of $H$:
  \begin{equation}\label{Eq:close2}
  \BoxCl(H) = \IIm{a^{(X)} \bp b^{(n - X)}}{X \ci n}.
  \end{equation}
Let $K$ be the element of $A \bp B$ given by the right hand side of
(\ref{Eq:clos}). We prove that $\BoxCl(H) = K$.

Let $X \in \Pow(n)$ and let $\fc \in \FD(n)$. If $X \in \fc$, then
$P_\fc(\vec a)\leq a^{(X)}$, while if $X \nin \fc$, then
$n - X \in \fc^*$, thus $P_{\fc^*}(\vec b)\leq b^{(n - X)}$. In both cases,
$P_\fc(\vec a)\ltp P_{\fc^*}(\vec b) \leq a^{(X)} \bp b^{(n - X)}$.
This proves that $K \ci\BoxCl(H)$.

Conversely, let $\vv<x,y> \in \BoxCl(H)$; we prove that $\vv<x,y> \in  K$.
If $x = 0_A$ or $y = 0_B$, then this is trivial, so suppose that both $x$
and $y$ are nonzero. Define
  \[
    \fc = \setm{X \ci n}{x \leq a^{(X)}} \ci \Pow(n).
  \]
It is trivial that $\fc$ is an upper subset of $\Pow(n)$. If $\fc =
\es$, then $n \nin \fc$, thus $x\nleq a^{(n)}$; but $\vv<x,y>  \in
a^{(n)} \bp b^{(\es)} = a^{(n)} \bp 0_B$, thus $y \leq 0_B$, a
contradiction. If $\fc = \Pow(n)$, then $\es  \in \fc$, thus $x \leq
a^{(\es)} = 0_A$, a contradiction.

Therefore, $\fc$ belongs to $\FD(n)$. By the definition of $\fc$, we have
$x \leq P_\fc(\vec a)$. Furthermore,
$n - X \nin \fc$, for all $X  \in \fc^*$,  which means that
$x\nleq a^{(n - X)}$. Since $\vv<x,y> \in a^{(n - X)} \bp b^{(X)}$,
the inequality $y \leq b^{(X)}$ holds.
This holds for all $X  \in \fc^*$, thus
$y \leq P_{\fc^*}(\vec b)$. Hence,
  \[
  \vv<x,y> \in P_\fc(\vec a) \ltp P_{\fc^*}(\vec b) \ci K,
  \]
which concludes the proof.
\end{proof}

Lemma~\ref{L:RepTens} implies two important purely arithmetical
formulas (see G.~A. Fraser \cite{Fras78} and~\cite{GrWe2}).

\begin{lemma}\label{L:IntersTens}
Let $A$ and $B$ be lattices with zero. Let $a_0$, $a_1 \in  A$
and $b_0$, $b_1 \in B$. Then
  \begin{align*}
(a_0 \ltp b_0) \ii (a_1 \ltp b_1)
  &= (a_0 \mm a_1) \ltp (b_0 \mm b_1),\\
(a_0 \ltp b_0) \jj (a_1 \ltp b_1) &= \\
(a_0 \ltp b_0) \uu(a_1 &\ltp b_1)
  \uu((a_0\jj a_1) \ltp(b_0 \mm b_1))
  \uu((a_0 \mm a_1) \ltp(b_0\jj b_1)).
  \end{align*}
  \end{lemma}

  \begin{proof}
  The first formula follows immediately from the definition of
$a \ltp b$.
The second formula  is a straightforward consequence
of Lemma~\ref{L:RepTens} (Formula (\ref{Eq:clos}), for $n = 2$).
  \end{proof}

If we further assume that either $a_0 \leq a_1$ and
$b_0 \geq b_1$, or $a_0 \geq a_1$ and $b_0 \leq b_1$, then
the second formula takes on the following simple form:
  \begin{equation}\label{Eq:TensId1}
(a_0 \ltp b_0)\jj(a_1 \ltp b_1) = (a_0 \ltp b_0) \uu(a_1 \ltp b_1).
  \end{equation}

\section{Semilattice tensor product and lattice tensor
product of lattices with zero}\label{S:stp}

For lattices $A$ and $B$ with zero, the \emph{extended \jz-semilattice
tensor product} $A\ootimes B$ is defined in \cite{GrWe2} as the set of all
\emph{bi-ideals} of $A\times B$ (see Definition~\ref{D:BotPureltp}(v)).
In particular, $A\ootimes B$ is an algebraic lattice.
The \jz-semilattice tensor product $A\otimes B$ is defined as the
\jz-semilattice of all compact elements of $A\ootimes B$.
The relationship between $A\otimes B$ (as in \cite{GLQu81},
\cite{GrWe2}, \cite{GrWe3} but not as in \cite{Fras78}) and the lattice
tensor product $A \ltp B$ is quite mysterious. Note that while
$A \otimes B$ may not be a lattice (see
  \cite{GrWe3} and \cite{GrWe4}), $A \ltp B$ is always a lattice.
Both $A\otimes B$ and $A \ltp B$ are \jz-semilattices.

\begin{corollary}\label{C:BoxRetr}
There exists a unique \jz-homomorphism $\rho$ from $A\otimes B$ to
$A\ltp B$ such that $\rho(a\otimes b)=a\ltp b$,
for all $\vv<a,b>\in A\times B$.
  \end{corollary}

Note that, in general, $A\ltp B$ is not a
join-subsemilattice of $A\otimes B$, even if $A\otimes B$ is a lattice.
\begin{proof}
We use the notation of \cite{GrWe2}. Every element $H$ of $A\ltp B$ is an
element of $A\bp B$, thus, by Lemma~\ref{L:biideal}, $H$ is a bi-ideal of
$A\times B$. Furthermore, by Lemma~\ref{L:withzero}, $H$ is a finite union
of pure lattice tensors. It follows that $H$ is a \emph{compact} element of
$A\ootimes B$, that is, an element of $A\otimes B$. Therefore,
$A \ltp B \ci A\otimes B$.

Let $a \in  A$ and let $b_0$, $b_1  \in  B$.
Since every element $H$ of $A\ltp B$ is a bi-ideal of $A\times B$,
  \[
     \vv<a,b_0\jj b_1> \in  H\quad \text{if{}f}\quad
      \vv<a,b_0>,\,\vv<a, b_1> \in H,
  \]
  from which it follows easily that
$a \ltp (b_0\jj b_1) = (a \ltp b_0)\jj(a \ltp b_1)$.
Furthermore, $a \ltp 0_B = 0_{A \ltp B}$.
By symmetry, it follows that the map from $A \times B$ to
$A \ltp B$ that sends every $\vv<a,b>$ to $a \ltp b$ is a \jz-bimorphism,
as defined in \cite{GrWe2}.
By the universal property of the tensor product, there exists a unique
\jz-homomorphism
$\gr\colon A\otimes B\to A \ltp B$ such that $\gr(a\otimes b) = a \ltp b$,
for any $\vv<a,b> \in  A \times B$. Therefore, $\gr$ is as desired.
  \end{proof}

\begin{proposition}\label{P:DLatTens}
Let $A$ and $B$ be lattices with zero. If either $A$ or $B$ is
distributive, then the semilattice tensor product and the lattice tensor
product of $A$ and $B$ coincide:
  \[
  A\otimes B = A \ltp B.
  \]
  \end{proposition}

Proposition~\ref{P:DLatTens} has an analogue for complete lattices, see, for
example, Corollary~5 in \cite{Wille}.

\begin{proof}
Without loss of generality, we can assume that $A$ is a distributive
lattice. Since $A \ltp B \ci A\otimes B$ always holds, we only have to
prove the converse. Let $H \in A\otimes B$; so there exists a decomposition
of the form
  \[
    H = \JJm{a_i \otimes b_i}{i < n}\qq(\text{computed in }A\otimes B),
  \]
  where $n$ is a positive integer and $\vv<a_i,b_i>
\in  A \times B$, for all $i < n$. Let $K$ be the corresponding element of
$A \ltp B$, that is,
  \[
     K = \JJm{a_i \ltp b_i}{i<n}\qq(\text{computed in }A\ltp B).
  \]
  We prove that $H = K$. Obviously, $H \ci K$. To prove the
converse, by Lemma~\ref{L:RepTens}, it suffices to prove that
  \[
  P_{\fc}(a_0,\ldots,a_{n-1}) \ltp P_{\fc^*}(b_0,\ldots,b_{n-1}) \ci H.
  \]
holds, for all $\fc \in \FD(n)$.

By Lemma 3.3 of \cite{GrWe3}, and Theorem 1 of \cite{Fras78}, it suffices to
prove that there exists a lattice polynomial $P$ such that
  \begin{equation}\label{Eq:Cruc}
    P_{\fc}(\vec a) \leq P(\vec a)\quad \text{and}\quad
     P_{\fc^*}(\vec b) \leq P^\dd(\vec b),
  \end{equation}
  where $P^\dd$ denotes the \emph{dual polynomial} of $P$.

We put $P=P_{\fc^*}^\dd$. Then $P^\dd=P_{\fc^*}$, thus
$P^\dd(\vec b)=P_{\fc^*}(\vec b)$. Since $A$ is distributive, it is easy
to verify that $P(\vec a)=P_\fc(\vec a)$. (Note that $P=P_\fc$ does not
hold in general; however, $P\leq P_\fc$.)
\end{proof}

\begin{remark}
  In Corollary 4.3 of \cite{GrWe3}, we proved that for all lattices $A$ and
$B$ with zero, if either $A$ or $B$ is distributive, then $A \otimes B$ is
a lattice.
  \end{remark}

\begin{example}\label{E:M3N5}
Denote by $M_3 = \set{0,p,q,r,1}$
and ${N}_5 = \set{0,a,b,c,1}$ (with $a>c$) the diamond and the
pentagon, respectively. We shall prove that
  \[
M_3 \ltp M_3 \ne  M_3\otimes M_3
\quad \text{and}\quad{N}_5 \ltp{N}_5 \ne {N}_5\otimes{N}_5.
  \]

Let $L$ be a finite lattice. We have seen in
  \cite{GrWe3} that there are natural isomorphisms
$\ga\colon M_3\otimes L\to M_3[L]$ and
$\ga'\colon {N}_5\otimes L\to {N}_5[L]$, where
$M_3[L]$ and ${N}_5[L]$ are the lattices defined by
  \begin{align}
M_3[L]& = \setm{\vv<x,y,z>  \in  L^3}{x \mm y = x \mm z = y \mm z},
\label{Eq:M3[L]}\\
{N}_5[L]& = \setm{\vv<x,y,z>  \in  L^3}{y \mm z \leq x \leq z}.
\label{Eq:N5[L]}
  \end{align}

The isomorphisms $\ga$ and $\ga'$ above are defined,
respectively, by the formulas
  \begin{align*}
    \ga(p \otimes x) &= \vv<x,0,0>, &\ga(q\otimes x) &= \vv<0,x,0>,
  &\ga(r\otimes x) &= \vv<0,0,x>;\\
    \ga'(a\otimes x) &= \vv<x,0,x>, &\ga'(b\otimes x) &= \vv<0,x,0>,
  &\ga'(c\otimes x) &= \vv<0,0,x>.
  \end{align*}

Define $M_3\vv<L>$ (resp., ${N}_5\vv<L>$) to be the image
of $M_3 \ltp L$ (resp., ${N}_5 \ltp L$) under $\ga$ (resp.,
$\ga'$).

Define the polynomials $\hat{x}$, $\hat{y}$, and $\hat{z}$
by $\hat{x} = y\jj z$, $\hat{y} = x\jj z$, and $\hat{z} = x\jj y$.
It is easy, though somewhat tedious, to compute that
  \begin{align}
    M_3\vv<L>& = \setm{\vv<x,y,z> \in  L^3}{x = \hat{y} \mm \hat{z},\
       y = \hat{x} \mm \hat{z},\ z = \hat{x} \mm \hat{y}},
      \label{Eq:M3<L>}\\
    {N}_5\vv<L>& = \setm{\vv<x,y,z> \in  L^3}{x = z \mm (x \jj y)}.
\label{Eq:N5<L>}
  \end{align}

\emph{In particular, $M_3\vv<L>$ has the same meaning here as
in \cite{GrWe1}}.

Thus it suffices to prove that $M_3\vv<M_3> \ne  M_3[M_3]$ and that
${N}_5\vv<{N}_5> \ne {N}_5[{N}_5]$. But it is easy to verify that
  \begin{align*}
  \vv<p,q,r>& \in M_3[M_3] - M_3\vv<M_3>,\\
  \vv<c,b,a>& \in {N}_5[{N}_5] - {N}_5\vv<{N}_5>.
  \end{align*}

By using \eqref{Eq:M3[L]} and \eqref{Eq:M3<L>}, it is also easy to see that
  \[
    M_3 \ltp {N}_5 = M_3\otimes {N}_5.
  \]
  \end{example}

\section{Lattice bimorphisms}

We shall see in this section one more reason to call the $A \ltp B$
construction the lattice tensor product.

\begin{definition}\label{D:bimorphism}
Let $A$, $B$, and $C$ be lattices with zero.
A \emph{$\set{0}$-lattice bimorphism} from $A \times B$ to $C$ is a map
$f\colon A \times B \to C$ such that
\begin{enumerate}
\item For all $\vv<a,b>  \in  A \times B$,
  \[
    f(\vv<a,0>) = f(\vv<0,b>) = 0.
  \]

\item For all $a_0$, $a_1  \in  A$ and all $b  \in  B$,
  \[
f(\vv<a_0\jj a_1,b>) = f(\vv<a_0,b>)\jj f(\vv<a_1,b>).
  \]

\item For all $a  \in  A$ and all $b_0$, $b_1  \in  B$,
  \[
f(\vv<a,b_0\jj b_1>) = f(\vv<a,b_0>)\jj f(\vv<a,b_1>).
  \]
\item For every positive integer $n$, all
$a_0$, \dots, $a_{n-1}$ in $A$, all $b_0$, \dots,
$b_{n-1}$ in $B$, and all $\fc  \in \FD(n)$,
  \[
f(\vv<P_\fc(a_0,\ldots,a_{n-1}),
P_{\fc^*}(b_0,\ldots,b_{n-1})>) \leq
\JJm{f(\vv<a_i,b_i>)}{i<n}.
  \]
\end{enumerate}
  \end{definition}

  Conditions (i)--(iii) define \jz-bimorphisms, see \cite{GrWe2}.
Condition (iv) is quite different, because it involves the meet structure
of $A$ and $B$ as well as the join structure.

\begin{proposition}\label{P:TensUniv}
Let $A$ and $B$ be lattices with zero. Consider the map
$ \ltp\colon A \times B\to A \ltp B$ defined by
$\vv<a, b>\mapsto a \ltp b$. Then $ \ltp$ is a universal
  $\set{0}$-lattice bimorphism, that is, for every lattice $C$ with zero
and every $\set{0}$-lattice bimorphism
$f \colon A \times B\to C$, there exists a unique
\jz-homomorphism $g\colon A \ltp B\to C$ such that
$g(a \ltp b) = f(\vv<a,b>)$, for all $a  \in  A$ and $b  \in  B$.
\end{proposition}

\begin{proof}
By Lemma~\ref{L:withzero}, the
elements of the form $a \ltp b$, where $\vv<a,b> \in A \times B$,
generate $A \ltp B$ as a \jz-semilattice. The uniqueness of $g$
follows immediately.

To prove the existence statement, it suffices to prove that
for every positive integer $n$, all $a$, $a_0$, \dots,
$a_{n-1}$ in $A$, and all $b$, $b_0$, \dots, $b_{n-1}$ in $B$,
  \begin{equation}\label{Eq:hyp}
  a \ltp b \leq\JJm{a_i \ltp b_i}{i<n}
  \end{equation}
implies that
  \begin{equation}\label{Eq:conc}
  f(\vv<a,b>) \leq\JJm{f(\vv<a_i,b_i>)}{i<n}.
  \end{equation}
The conclusion (\ref{Eq:conc}) is trivial if $a = 0_A$ or
$b = 0_B$, so suppose that both $a$ and $b$ are nonzero.
By Lemma~\ref{L:RepTens}, (\ref{Eq:hyp}) is equivalent to the
existence of an element $\fc$ of $\FD(n)$ such that
  \[
  a \leq P_\fc(a_0,\ldots,a_{n-1})\quad \text{and}\quad
  b \leq P_{\fc^*}(b_0,\ldots,b_{n-1}).
  \]
Since $f$ is a \jz-bimorphism, it is isotone, thus
   \begin{align*}
  f(\vv<a,b>)& \leq f(\vv<P_\fc(a_0,\ldots,a_{n-1}),
  P_{\fc^*}(b_0,\ldots,b_{n-1})>)\\
  & \leq\JJm{f(\vv<a_i,b_i>)}{i<n},
  \end{align*}
because $f$ is a $\set{0}$-lattice bimorphism, which completes the proof.
\end{proof}

This shows that $\ltp$ defines, in fact, a \emph{bifunctor} on $\E L_0$. A
useful direct description of the effect of this functor on morphisms in
$\E L_0$ is given by the following result.

\begin{proposition}\label{P:CanHom}
Let $A$, $A'$, $B$, $B'$ be objects in $\E L_0$ and let $f\colon A\to A'$
and $g\colon B\to B'$ be morphisms in $\E L_0$.  Then
  $(f \ltp g)(X)$ is given by the following formula, for all $X  \in  A
\ltp B$:
  \begin{equation}\label{Eq:TensMap}
    (f \ltp g)(X)=\UUm{f(x)\ltp g(y)}{\vv<x,y>\in X}.
  \end{equation}
  \end{proposition}

\begin{proof}
Let $h$ be the map defined on the powerset of $A \times B$ by
the formula (\ref{Eq:TensMap}); denote by $h'$ the
restriction of $h$ to $A \ltp B$. It suffices to prove that
$h' = f \ltp g$.

Since $f$ and $g$ are morphisms in $\E L_0$,
  \begin{equation}\label{Eq:hPureTens}
h(a \ltp b) = f(a) \ltp g(b).
  \end{equation}
holds, for all $\vv<a,b> \in  A \times B$. Let $X$ be an arbitrary element
of $A \ltp B$. There exists a decomposition of $X$ of the form
  \[
  X = \JJm{a_i \ltp b_i}{i<n},
  \]
where $n$ is a positive integer and $\vv<a_i,b_i> \in  A \times B$, for
all $i$. By Lemma~\ref{L:RepTens},
  \begin{equation}\label{Eq:RepH}
  X = \UUm{P_{\fc}(a_0,\ldots,a_{n-1}) \ltp
  P_{\fc^*}(b_0,\ldots,b_{n-1})}{\fc \in \FD(n)}.
  \end{equation}
But, by definition, $h$ is a join-homomorphism from
$\Pow(A \times B)$ to $\Pow(A' \times B')$. Therefore, it
follows from \eqref{Eq:hPureTens}, \eqref{Eq:RepH},
and the fact that $f$ and $g$ are morphisms in $\E L_0$
that
  \[
   h(X) =
  \UUm{P_{\fc}(f(a_0),\ldots,f(a_{n-1})) \ltp
      P_{\fc^*}(g(b_0),\ldots,g(b_{n-1}))}{\fc \in \FD(n)},
  \]
thus, again by Lemma~\ref{L:RepTens},
  \[
  h(X) = \JJm{f(a_i) \ltp g(b_i)}{i<n}.
  \]
  We conclude that $h' = f \ltp g$.
\end{proof}

As an immediate corollary, every object of $\E L_0$ is \emph{flat} with
respect to the lattice tensor product bifunctor $\ltp$:

\begin{proposition}\label{P:Flat}
In the context of Proposition~\tup{\ref{P:CanHom}}, if both $f$ and $g$ are
lattice embeddings, then so is $f \ltp g$.
  \end{proposition}

Another fact worth mentioning is that $f \ltp g$ is a restriction of
$f\otimes g$:

\begin{corollary}\label{C:fgotltp}
In the context of Proposition~\tup{\ref{P:CanHom}}, $f \ltp g$ is the
restriction from $A \ltp B$ to $A' \ltp B'$ of the map $f\otimes g\colon
A\otimes B\to A'\otimes B'$.
  \end{corollary}

\begin{proof}
This is an immediate consequence of Proposition 3.4 of \cite{GrWe2}.
\end{proof}

\section{$A \ltp B$ as a capped sub-tensor product}

In \cite{GrWe2}, we introduced the following definition:

\begin{definition}\label{D:capped}
  Let $A$ and $B$ be lattices with zero. A \emph{sub-tensor
product} of $A$ and $B$ is a subset $C$ of the semilattice tensor
product $A \otimes B$ satisfying the following conditions:

  \begin{enumerate}
\item $C$ is closed under finite intersection.

\item $C$ is a lattice under containment.

\item For all $a_0$, $a_1 \in A$ and all $b_0$, $b_1 \in B$, if
either $a_0 \leq a_1$ and $b_0 \geq b_1$, or
$a_0 \geq a_1$ and $b_0 \leq b_1$, then the hereditary set
  \[
    (a_0\otimes b_0) \uu(a_1\otimes b_1)\q\q (\emph{mixed tensor})
  \]
  belongs to $C$.
  \end{enumerate}

A \emph{capped sub-tensor product} of $A$ and $B$ is a sub-tensor product
of $A$ and $B$ satisfying the following additional condition:

  \begin{enumerate}
\item[(iv)] Every element of $C$ is a finite \emph{union} of pure tensors.
  \end{enumerate}
  \end{definition}

It is an open problem whether every sub-tensor product is capped, see
Problem~2 in \cite{GrWe2}.

$A \ltp B$ is an example of a capped sub-tensor product:

\begin{theorem}\label{T:ltpTensLatt}
Let $A$ and $B$ be lattices with zero. Then $A \ltp B$ is a capped
sub-tensor product of $A$ and $B$. Furthermore, it is the smallest
(with respect to containment) sub-tensor product of $A$ and $B$.
  \end{theorem}

\begin{proof}
By Proposition~\ref{P:AltpBlatt}, $A\ltp B$ is an ideal of $A\bp B$. Since
$A\bp B$ is a lattice under containment (Proposition~\ref{P:bplattice}),
closed under finite intersection,
$A \ltp B$ satisfies (i) and~(ii).
  Furthermore, (iii) follows immediately from the
particular case \eqref{Eq:TensId1} of Lemma~\ref{L:IntersTens}.
  Finally, (iv) follows from Lemma~\ref{L:withzero}.

Now let $C$ be a sub-tensor product of $A$ and $B$; we prove that $C$
contains $A\ltp B$. So let $H\in A\ltp B$. Then $H$ belongs to $A\bp B$,
thus $H$ can be written in the following form:
  \[
  H=\IIm{a_i\bp b_i}{i<n},
  \]
where $n$ is a positive integer and $\vv<a_i,b_i>\in A\times B$.
Furthermore, $H$ is confined, thus there exists
$\vv<a,b>\in A\times B$ such that $H\ci a\ltp b$. Hence,
  \[
  H=\IIm{(a_i\bp b_i)\ii(a\ltp b)}{i<n}.
  \]
However, for all $i<n$, it is easy to compute that
  \[
  (a_i\bp b_i)\ii(a\ltp b)=
  ((a\mm a_i)\ltp b)\uu(a\ltp(b\mm b_i)),
  \]
which is a mixed tensor. Therefore, by the definition of a sub-tensor
product, $H$ belongs to $C$.
  \end{proof}

We can then use Theorem 2 of \cite{GrWe2} to deduce the
following result:

\begin{theorem}\label{T:AotBCon}
Let $A$ and $B$ be lattices with zero. Then there exists a
unique isomorphism $\gm$ from $\Conc A\otimes\Conc B$ onto
$\Conc(A \ltp B)$ such that, for all $a_0 \leq a_1$ in $A$ and
all $b_0 \leq b_1$ in $B$, the following equality holds:
  \begin{equation*}
    \gm(\gQ_A(a_0,a_1)\otimes\gQ_B(b_0,b_1)) =
    \gQ_{A \ltp B}((a_0 \ltp b_1)\jj(a_1 \ltp b_0),a_1 \ltp b_1).
  \end{equation*}
\end{theorem}

Theorem A follows immediately.

\section{$\set{1}$-sensitive homomorphisms; the box product bifunctor}

The box product operation, $\bp$, is not a bifunctor from the category of
lattices with lattice homomorphisms to itself. However, we will see that
considering only the following general type of homomorphism will overcome
this difficulty.

\begin{definition}\label{D:1sep}
Let $A$, $B$ be lattices, let $f\colon A\to B$ be a lattice
homomorphism. We will say that $f$ is \emph{$\set{1}$-sensitive}, if
$1_A$ exists if and only if $1_B$ exists, and if they both exist then
$f(1_A) = 1_B$.
  \end{definition}

Note that if $f\colon A\to B$ is a lattice homomorphism and
neither $1_A$ nor $1_B$ exists, then $f$ is $\set{1}$-sensitive.

It is clear that lattices and $\set{1}$-sensitive maps form a
subcategory of the category of all lattices and lattice
homomorphisms.

\begin{proposition}\label{P:1PresFunc}
Let $A$, $A'$, $B$, and $B'$ be lattices, let $f\colon A\to A'$
and $g\colon B\to B'$ be $\set{1}$-sensitive lattice homomorphisms.
Then there exists a unique map
$h$ from $A \bp B$ to $A' \bp B'$ such that
  \begin{equation}\label{Eq:DeffbpC}
  h\left(\IIm{a_i \bp b_i}{i<n}\right) =
  \IIm{f(a_i) \bp g(b_i)}{i<n}.
  \end{equation}
  holds, for every positive integer $n$ and all $a_i \in  A$,
$b_i \in B$ ($i<n$). Furthermore, $h$ is a $\set{1}$-sensitive lattice
homomorphism.
  \end{proposition}

\begin{proof}
The uniqueness statement is trivial. To prove existence of a
map $h$ satisfying (\ref{Eq:DeffbpC}), it is sufficient to
prove that
  \begin{equation}\label{Eq:HypInc}
     \IIm{a_i \bp b_i}{i<n} \ci a \bp b
  \end{equation}
  implies that
  \begin{equation}\label{Eq:ConcInc}
     \IIm{f(a_i) \bp g(b_i)}{i<n} \ci f(a) \bp g(b),
  \end{equation}
  for all $n>0$ and all $a$, $a_i \in  A$, $b$, $b_i \in  B$ ($i<n$).
Now \eqref{Eq:HypInc} is equivalent to the following condition:
  \begin{gather*}
  a = 1_A \text{\q or\q }b = 1_B\\
   \intertext{or}
  a_{(n)} \leq a \text{\q and\q }b_{(n)} \leq b \text{\q and\q }
  (\forall X  \in \tup{P}^*(n))
  (a_{(X)} \leq a \text{ or }b_{(n -  X)} \leq b).
  \end{gather*}

\medskip

\noindent Since $f$ and $g$ are $\set{1}$-sensitive lattice homomorphisms,
this implies the condition:
  \[
  f(a) = 1_{A'} \text{\q or\q }g(b) = 1_{B'}
\]
  or
  \begin{align*}
    a'_{(n)} \leq f(a) &\text{\q and\q }b'_{(n)} \leq g(b)\\
    &\text{\q and\q } (\forall X  \in \tup{P}^*(n))
      (a'_{(X)} \leq f(a) \text{ or }b'_{(n -  X)} \leq g(b)).
  \end{align*}
\noindent which, in turn, is equivalent to \eqref{Eq:ConcInc}.

We now verify that $h$ is a lattice homomorphism. It is obvious that $h$
is a meet homomorphism. The fact that $h$ is a join homomorphism follows
immediately from Lemma~\ref{L:CompJoin}.

Since both $f$ and $g$ are $\set{1}$-sensitive, $1_A$ exists if and only
if $1_{A'}$ exists, and $1_B$ exists if and only if $1_{B'}$ exists. By
Remark~\ref{R:1AboxB}, $1_{A \bp B}$ exists if and only if $1_{A' \bp B'}$
exists. Suppose now that $1_{A \bp B}$ and $1_{A' \bp B'}$
exist. Without loss of generality, $1_A$~exists. Since $f$ is
$\set{1}$-sensitive, $1_{A'}$ exists and $f(1_A) = 1_{A'}$, thus
  \[
    h(1_{A \bp B}) = h(1_A \bp b) = f(1_A) \bp g(b) = 1_{A'} \bp g(b) =
       A' \times B',
  \]
  for all $b \in B$, and so $1_{A' \bp B'}$ exists. Therefore, $h$ is
$\set{1}$-sensitive.
\end{proof}

We shall denote by $f \bp g$ the $\set{1}$-sensitive lattice homomorphism
$h$ of Proposition~\ref{P:1PresFunc}.

\begin{remark}
In the proof of Proposition~\ref{P:1PresFunc}, in order to prove the
existence of a lattice homomorphism $h$ satisfying (\ref{Eq:DeffbpC}), we
require only a weaker assumption on $f$ and $g$: namely, if $1_A$ exists,
then $1_{A'}$ exists and $f(1_A) = 1_{A'}$. However, we shall require
later the stronger definition of a $\set{1}$-sensitive map for direct
limits (see Proposition~\ref{P:AltpPresLim}).
  \end{remark}

The following consequence of
Proposition~\ref{P:1PresFunc} is immediate:

\begin{corollary}\label{C:1PresFunc}
The mappings $\vv<A,B>\mapsto A \bp B$,
$\vv<f,g>\mapsto f \bp g$ define a bifunctor from the category
of lattices and $\set{1}$-sensitive lattice homomorphisms to
itself.
  \end{corollary}

The following corollary will be of special importance:

\begin{corollary}\label{C:bpltpFun}
Let $A$, $B$, and $C$ be lattices, with $A$ bounded, and let $f\colon B\to
C$ be a $\set{1}$-sensitive lattice homomorphism. Then the image of $A
\ltp B$ under $\id_A \bp f$ is contained in $A \ltp C$.
  \end{corollary}

\begin{proof}
Put $g = \id_A \bp f$. We prove that $g(H) \in  A \ltp C$, for all $H
\in  A \ltp B$. By the definition of $A \ltp B$, one can write $H$ in the
form
  \[
H = \IIm{a_i \bp b_i}{i<n},
  \]
where $n>0$, $a_i  \in  A$, $b_i  \in  B$ (for all $i<n$), and
$\MMm{a_i}{i<n} = 0_A$. Therefore, we obtain that
  \[
g(H) = \IIm{a_i \bp f(b_i)}{i<n}.
  \]
Since $\MMm{a_i}{i<n} = 0_A$, we conclude, by Lemma~\ref{L:Abdd},
that $g(H)$ belongs to $A \ltp C$.
  \end{proof}

In the context of Corollary~\ref{C:bpltpFun}, we will write
$\id_A \ltp f$ for the restriction of $\id_A \bp f$ from $A \ltp
B$ to $A \ltp C$. Similarly, we define $f \ltp\id_C$, if $C$
is a bounded lattice and $f\colon A\to B$ is a $\set{1}$-sensitive
lattice homomorphism.

\section{The functor $A \ltp\ghost$, for $A$ bounded}
\label{S:DirLim}

In this section, we investigate box products of lattices where one of the
factors is bounded.

\begin{proposition}\label{P:AltpPresLim}
Let $A$ be a bounded lattice. Let $\vv<I, \leq>$ be a
directed set, let $B$, $B_i$ ($i  \in  I$) be lattices such that, for
appropriate $\set{1}$-sensitive transition maps $f_{ij}\colon B_i\to B_j$
(for $i \leq j$) and $f_i\colon B_i\to B$, we have
  \[
  B = \varinjlim_iB_i.
  \]
Then, with the transition maps $g_{ij} = \id_A \ltp f_{ij}$ and
$g_i = \id_A \ltp f_i$, we have
  \[
  A \ltp B = \varinjlim_iA \ltp B_i.
  \]
\end{proposition}

\begin{proof}
It suffices to prove that for all $i  \in  I$ and for all $H$, $K  \in  A
\ltp B_i$, $g_i(H) \ci g_i(K)$ implies that there exists $j \geq i$ in $I$
such that $g_{ij}(H) \ci g_{ij}(K)$. Write $H$ and $K$ as
  \begin{align*}
    H& = \IIm{a_k \bp b_k}{k<m},
  \end{align*}
  where $m>0$ and $\MMm{a_k}{k<m} = 0_A$, and
  \begin{align*}
     K& = \IIm{c_l \bp d_l}{l<n},
  \end{align*}
  where $n>0$  and $\MMm{c_l}{l<n} = 0_A$.

The assumption $g_i(H) \ci g_i(K)$ means that
  \begin{equation}\label{Eq:akbkcldl1}
     \IIm{a_k \bp f_i(b_k)}{k<m} \ci c_l \bp f_i(d_l),
  \end{equation}
  holds, for all $l<n$. Since $I$ is directed, it suffices to prove that for
all $l<n$ there exists $j \geq i$ in $I$ such that
  \begin{equation}\label{Eq:akbkcldl2}
     \IIm{a_k \bp f_{ij}(b_k)}{k<m} \ci c_l \bp f_{ij}(d_l).
  \end{equation}
If $c_l = 1_A$, then this is trivial (take $j = i$). If $f_i(d_l) = 1_B$,
then, since $f_i$ is $\set{1}$-sensitive, $1_{B_i}$ exists and
$f_i(1_{B_i}) = 1_B$. It follows that $f_i(1_{B_i}) = f_i(d_l)$, thus
there exists $j \geq i$ in $I$ such that $f_{ij}(1_{B_i}) = f_{ij}(d_l)$.
Since $f_{ij}$ is $\set{1}$-sensitive, it follows that $f_{ij}(d_l) =
1_{B_j}$; (\ref{Eq:akbkcldl2}) follows. Suppose now that $c_l$ is not the
largest element of~$A$, and that $f_i(d_l)$ is not the largest element of
$B$. Then (\ref{Eq:akbkcldl1}) means that, for all $X  \in \tup{P}^*(m)$,
either $a_{(X)} \leq c_l$ or $f_i(b_{(m -  X)}) \leq f_i(d_l)$. Since $B =
\varinjlim_jB_j$, we obtain that there exists $j \geq i$ in $I$ such that
the conditions above hold with $f_{ij}$ instead of $f_i$. Then
(\ref{Eq:akbkcldl2}) follows; whence $g_{ij}(H) \ci g_{ij}(K)$.
  \end{proof}

For every lattice $L$ with zero, denote by $\gl_L$ the canonical
isomorphism from $\Conc A\otimes\Conc L$ onto $\Conc(A \ltp L)$. Define
the functors, $\gF$ and $\gY$, from lattices and $\set{1}$-sensitive
homomorphisms to semilattices with zero and \jz-homomorphisms, by
\begin{align*}
    \gF(L) &= \Conc A\otimes\Conc L,\\
    \gY(L) &= \Conc(A \ltp L),
\end{align*}
  extended to morphisms in the natural way.

\begin{lemma}\label{L:NatEqltp}
Let $A$ be a bounded lattice. The correspondence
$L\mapsto\gl_L$ defines a natural transformation from the functor $\gF$ to
the functor $\gY$ on the subcategory of lattices with zero.
  \end{lemma}

\begin{proof}
This amounts to verifying, for $f\colon B\to C$ a
$\set{1}$-sensitive homomorphism of lattices with zero, that the
following diagram
  \[
    \begin{CD}
      \gY(B) @> \gY(f) >> \gY(C)\\
      @A \gl_B AA @AA \gl_C A\\
      \gF(B) @> \gF(f) >> \gF(C)
     \end{CD}
  \]
is commutative. It suffices to prove that every congruence of
the form
  \[
    \gQ = \gQ_A(a_0,a_1) \otimes\gQ_B(b_0,b_1),
  \]
  where $a_0 \leq a_1$ in $A$ and $b_0 \leq b_1$ in $B$, has the same
image under the maps $\gl_C \circ \gF(f)$ and $\gY(f) \circ \gl_B$. We
compute:
  \begin{align*}
    \gY(f) \circ\gl_B(\gQ)& = \gY(f)\bigl(\gQ_{A \ltp B}
   ((a_0 \ltp b_1)\jj(a_1 \ltp b_0),a_1 \ltp b_1) \bigr)\\
     & = \gQ_{A \ltp C} ((a_0 \ltp f(b_1))\jj(a_1 \ltp f(b_0)),a_1 \ltp
        f(b_1)),\\
   \intertext{while}
   \gl_C \circ\gF(f)(\gQ)& = \gl_C
    \bigl(
      \gQ_A(a_0,a_1)\otimes\gQ_C(f(b_0),f(b_1))
    \bigr)\\
    & = \gQ_{A \ltp C} ((a_0 \ltp f(b_1))\jj(a_1 \ltp f(b_0)),a_1 \ltp
         f(b_1)),
  \end{align*}
which concludes the proof.
\end{proof}

We can now deduce the following extension of
Theorem~\ref{T:AotBCon}:

\begin{theorem}\label{T:AotBCon2}
Let $A$ and $B$ be lattices, with $A$ bounded. Then there
exists a unique isomorphism $\gm$ from
$\Conc A\otimes\Conc B$ onto $\Conc(A \ltp B)$ such that
  \[
     \gm(\gQ_A(a_0,a_1)\otimes\gQ_B(b_0,b_1)) =
     \gQ((a_0 \bp b_0) \ii (0_A \bp b_1),(a_1 \bp b_0) \ii (0_A \bp b_1)).
  \]
  holds, for all $a_0 \leq a_1$ in $A$ and $b_0 \leq b_1$ in $B$.
\end{theorem}

Note that, indeed, both elements $(a_0 \bp b_0) \ii (0_A \bp b_1)$ and
$(a_1 \bp b_0) \ii (0_A \bp b_1)$ belong to $A \ltp B$.

\begin{proof}
  The uniqueness of $\gm$ is obvious. To prove the existence, we represent
$B$ as the direct limit of all its sublattices $B_b = [b)$, for $b  \in
B$; the index set is the partially ordered set dual of $B$, the transition
maps are all the inclusion maps. They are obviously $\set{1}$-sensitive.
Therefore, the following isomorphisms hold, with the canonical transition
maps:
   \begin{align*}
     \Conc(A \ltp B)&\cong\varinjlim_b\Conc(A \ltp B_b)\\
   \intertext{(by Proposition~\ref{P:AltpPresLim} and the fact
      that the functor $\Conc$ preserves direct limits)}
    &\cong\varinjlim_b\Conc A\otimes\Conc B_b\\
    \intertext{(by Lemma~\ref{L:NatEqltp})}
    &\cong\Conc A\otimes\Conc B
  \end{align*}
(because the functors $\Conc$ and
     $\Conc A\otimes\ghost$ preserve direct limits).
Denote by $\gm\colon\Conc A\otimes\Conc B\to\Conc(A \ltp B)$
the isomorphism thus obtained. We compute the effect of $\gm$
on $\gQ = \gQ_A(a_0,a_1)\otimes\gQ_B(b_0,b_1)$, with
$a_0 \leq a_1$ in $A$ and $b_0 \leq b_1$ in $B$.

Put $b = b_0$, and $\gQ' = \gQ_A(a_0,a_1)\otimes\gQ_{B_b}(b_0,b_1)$. Keep
the notations $\gF$, $\gY$ for the two functors defined above (with
parameter $A$), and $L\mapsto\gl_L$ for the natural transformation from
$\gF$ to $\gY$. Put $g_b = \id_A \ltp f_b$. Then we compute
   \begin{align*}
  \gm(\gQ)& = \gm \circ\gF(f_b)(\gQ')\\
  & = \gY(f_b) \circ\gl_{B_b}(\gQ')\\
  & = \gY(f_b)(
  \gQ_{A \ltp B_b}((a_0 \ltp b_1)\jj(a_1 \ltp b_0),a_1 \ltp b_1)
  )\\
  & = \gQ_{A \ltp B}(
  g_b((a_0 \ltp b_1)\jj(a_1 \ltp b_0)),g_b(a_1 \ltp b_1)).
  \end{align*}
It is not difficult to compute that, in $A \ltp B_b$, we have
   \begin{align*}
  (a_0 \ltp b_1)\jj(a_1 \ltp b_0) = (a_0 \bp b_0) \ii (0_A \bp b_1),\\
   \intertext{while}
  a_1 \ltp b_1 = (a_1 \bp b_0) \ii (0_A \bp b_1).
  \end{align*}
The conclusion follows.
\end{proof}

Theorem B follows immediately. However, we could not find a construction
proving that the tensor product of two representable join semilattices with
zero is again representable. See also Problems \ref{Pb:Repr},
\ref{Pb:CRepr} and \ref{Pb:CTRepr}.

It is easy to deduce the following far reaching generalization of the main
result of \cite{GrWe1}:

\begin{corollary}\label{C:ConPres}
Let $S$ and $L$ be lattices, with $S$ bounded and simple. Then $L$ admits
a congruence-preserving embedding into $S \ltp L$, defined by
$x\mapsto 0_S \bp x$.
  \end{corollary}

\section{Congruences on box product of lattices with unit}

A similar direct limit argument as the one used in Section~\ref{S:DirLim}
yields a result about congruences on box products of lattices with unit,
similar to Theorems \ref{T:AotBCon} and~\ref{T:AotBCon2}. However, there
is a much less painful way of obtaining this.

\begin{theorem}\label{T:AotBCon3}
Let $A$ and $B$ be lattices with unit. Then there exists a
unique isomorphism $\gm$ from $\Conc A\otimes\Conc B$ onto
$\Conc(A \bp B)$ such that
  \begin{equation*}
  \gm(\gQ_A(a_0,a_1)\otimes\gQ_B(b_0,b_1)) =
  \gQ_{A \bp B}(a_0 \bp b_0,(a_0 \bp b_1) \ii (a_1 \bp b_0)),
  \end{equation*}
for all $a_0 \leq a_1$ in $A$ and all $b_0 \leq b_1$ in $B$.
  \end{theorem}

\begin{proof}
The following isomorphisms hold:
  \begin{align*}
    \Conc A\otimes\Conc B &\cong\Conc A^\dd\otimes\Conc B^\dd\\
    &\cong\Conc(A^\dd \ltp B^\dd) &&(\text{by
        Theorem~\ref{T:AotBCon}})\\
    &\cong\Conc(A \bp B)^\dd &&(\text{by Corollary~\ref{C:BpLtpDual}})\\
    &\cong\Conc(A \bp B).
  \end{align*}
  Furthermore, the successive images of the tensor product of two
principal congruences $\gQ_A(a_0,a_1)$ and $\gQ_B(b_0,b_1)$
(with $a_0 \leq a_1$ and $b_0 \leq b_1$) under the
isomorphisms above are the following (see Remark~\ref{R:BpLtpDual}):
  \begin{align*}
    \gQ_A(a_0,a_1)\otimes\gQ_B(b_0,b_1) &\mapsto
    \gQ_{A^\dd}(a_1,a_0)\otimes\gQ_{B^\dd}(b_1,b_0)\\
      &\mapsto \gQ_{A^\dd \ltp B^\dd}((a_0 \ltp b_1)\jj(a_1 \ltp b_0),a_0
     \ltp b_0)\\
     &\mapsto \gQ_{(A \bp B)^\dd} ((a_0 \bp b_1)\jj(a_1 \bp b_0),a_0 \bp
       b_0)\\
     &\mapsto \gQ_{A \bp B}
       (a_0 \bp b_0,(a_0 \bp b_1) \ii (a_1 \bp b_0)),
  \end{align*}
  which proves the existence statement. The uniqueness is obvious.
\end{proof}

\section{Discussion}\label{S:discussion}

The various tensor products of lattices show an interesting formal
similarity among some of the results. These constructions:
\begin{enumerate}
\item  preserve distributivity (of lattices or of semilattices);
\item  can be characterized with maps from one lattice to the other;
\item  have an ``Isomorphism Theorem'' for their (compact) congruence
(semi) lattices.
\end{enumerate}

We refer to B. Ganter and R. Wille \cite{GaWi},
G. Gr\"atzer, H Lakser, and R.~W. Quackenbush \cite{GLQu81},
R.~W. Quackenbush \cite{rQ85}, G.~N. Raney \cite{Rane}, Z. Shmuely
\cite{Shmu}, R. Wille \cite{Wille}, and the authors' papers
\cite{GrWe1}--\cite{GrWe4}, for more information.

More interestingly, it seems that formally similar results for two different
types of tensor products do not seem to imply each other. For
example, consider the Isomorphism Theorem for compact
congruence semilattices of tensor products of lattices (Theorem 2 of
\cite{GrWe2}):
  \begin{equation}\label{Eq:IsoThmCapped}
  \Conc(A\otimes B)\iso\Conc A\otimes\Conc B,
  \end{equation}
provided that $A$ and $B$ are lattices with zero and $A\otimes B$ is a
lattice and the Isomorphism Theorem for complete
congruence lattices of doubly founded complete lattices (Theorem~18 in
\cite{Wille}):
  \begin{equation}\label{Eq:IsoThmCpl}
  \Coni(A\hotimes B)\iso\Coni A\hotimes\Coni B,
  \end{equation}
where $A\hotimes B$ is the complete tensor product introduced in R.~Wille
\cite{Wille}, and $\Coni K$ is the complete congruence lattice of a complete
lattice $K$.

Both results apply to finite lattices. For finite lattices $A$ and $B$,
Wille's Isomorphism Theorem is a special case of
Theorem~\ref{T:AotBCon}, which is similar, though not equivalent, to the
Isomorphism Theorem for tensor products of finite lattices in \cite{GLQu81}.
For infinite lattices $A$ and $B$, the two Isomorphism Theorems seem to have
nothing in common:  \eqref{Eq:IsoThmCapped} equates tensor
products of two \emph{distributive} \jz-semilattices, while
\eqref{Eq:IsoThmCpl} equates tensor products of arbitrary
complete lattices. It was proved in G. Gr\"atzer \cite{gG90} (see G.
Gr\"atzer and H. Lakser \cite{GrLa91} for the shortest proof and G.~
Gr\"atzer and E.~T. Schmidt \cite{GS95} for the strongest result)
that $\Coni A$ can be any complete lattice.

In general, the constructions of complete tensor products of complete
lattices are given as complete meet-semilattices, so, of course, they are
lattices. The situation is quite different for tensor product constructions
of (not necessarily complete) lattices, where the tensor product may
\emph{not} be a lattice, see \cite{GrWe3} and \cite{GrWe4}. So, in one sense,
Proposition~\ref{P:bplattice} lies at the core of the present paper.

This difficulty is paralleled in the characterization problems of congruence
lattices: while complete congruence lattices of complete lattices have been
characterized, see \cite{gG90}, the characterization problem of congruence
lattices of lattices is open, see G. Gr\"atzer and E.~T. Schmidt
\cite{GrSc} for a survey.

\section{Open problems}

\begin{problem}
  Denote by $\mbf{V}(L)$ the variety generated by a lattice $L$. Let $A$
and $B$ be lattices with zero.  Prove that $A \ltp B = A\otimes B$ if and
only if $\mbf{V}(A) \ii \mbf{V}(B)$ is a distributive variety.
  \end{problem}

See Example~\ref{E:M3N5} for some basic examples related to
this problem.

\begin{problem}\label{Pb:Repr}
Is every representable semilattice $\set{0}$-representable?
\end{problem}

It would follow, by Theorem A, that the tensor product of any two
representable distributive semilattices with zero is representable. On the
other hand, it is not even known whether there exists a nonrepresentable
distributive semilattice with zero.

However, the situation changes if we consider lattices \emph{with permutable
congruences}.  Let us say that a \jz-semilattice $D$ is
\emph{\textup{p}-representable} (resp.,
$\vv<\textup{p},\set{0}>$-representable), if there exists a lattice (resp.,
a lattice with zero) $L$ with permutable congruences such that $\Conc L\iso
D$. There are non \textup{p}-representable distributive \jz-semilattices,
see J.~T\r{u}ma and F. Wehrung \cite{TuWe}. Furthermore, the second author
of the present paper proved the following result: \begin{quote}
\em Let $A$ and $B$ be lattices with permutable congruences. If
$A\ltp B$ is defined, then $A\ltp B$ has permutable congruences.
\end{quote}

In particular, if $S$ and $T$ are $\vv<\textup{p},\set{0}>$-representable
\jz-semilattices, then $S\otimes T$ is
$\vv<\textup{p},\set{0}>$-representable. Hence a reasonable analogue of
Problem~\ref{Pb:Repr} for lattices with permutable congruences is the
following:

\begin{problem}\label{Pb:CRepr}
Is every \textup{p}-representable semilattice
$\vv<\textup{p},\set{0}>$-representable?
  \end{problem}

A problem more directly related to tensor products is the following:

\begin{problem}\label{Pb:CTRepr}
If $S$ and $T$ are \textup{p}-representable \jz-semilattices, is
$S\otimes T$ \textup{p}-rep\-re\-sent\-a\-ble?
\end{problem}

Any counterexample to Problem~\ref{Pb:CTRepr} must have either $S$ or
$T$ not $\vv<\textup{p},\set{0}>$-represent\-able and either $S$ or
$T$ must have at least $\aleph_2$ elements. Such a result would imply a
negative answer to Problem~\ref{Pb:CRepr}.

  \begin{problem}
  Are there other lattice tensor product constructions between $A \ltp B$ and
$A\otimes B$? For example, in view of Lemma~\ref{L:RepTens}, we could assert
that the $A \ltp B$ construction utilizes the structure of the free
distributive lattices. Are there analogues of $A \ltp B$ for other varieties
of lattices?
  \end{problem}

If $A$ and $B$ are lattices with zero, then $A \ltp B$
is the smallest capped sub-tensor product of $A$ and $B$ (see
Theorem~\ref{T:ltpTensLatt}). On the other hand, if $A\otimes B$ is a
capped tensor product, then $A\otimes B$ is the largest capped
sub-tensor
product of $A$ and $B$.

\begin{problem}
The tensor
product of two finite simple lattice is a larger
finite simple
lattice. In general, what are the ``ultimate
building blocks" of,
say, finite lattices, by using
elementary operations such as direct
product, ordinal sum,
and generalizations of the tensor
product?
\end{problem}

\begin{problem}
What can be said about
\emph{relative tensor products}, that
is, lattice-theoretical
analogues of the module-theoretical
construction $A\otimes_RB$? Does there exist such a
construction?
\end{problem}

\section*{Acknowledgment}
This work was partially completed
while the second author was visiting the University of
Manitoba. The excellent conditions provided by the Mathematics
Department, and, in particular, a quite lively seminar, were
greatly appreciated.

We wish to thank the  referee for his constructive suggestions.

\end{document}